\documentclass[12]{article}
\usepackage{fleqn}
\usepackage{graphicx}
\begin{document}
\Large
\begin{center}
{\bf A Classification of the Projective Lines over Small Rings}
\end{center}
\vspace*{.5cm}
\begin{center}
Metod Saniga,$^{1-3}$ Michel Planat,$^{2}$ Maurice R. Kibler$^{3}$\\
and Petr Pracna$^{2-4}$
\end{center}
\vspace*{.2cm} \normalsize
\begin{center}
$^{1}$Astronomical Institute, Slovak Academy of Sciences\\
SK-05960 Tatransk\' a Lomnica, Slovak Republic\\
(msaniga@astro.sk)

\vspace*{.2cm} $^{2}$Institut FEMTO-ST, CNRS, D\' epartement LPMO,
32 Avenue de
l'Observatoire\\ F-25044 Besan\c con Cedex, France\\
(planat@lpmo.edu)

\vspace*{.2cm} $^{3}$Institut de Physique Nucl\' eaire de Lyon,
IN2P3-CNRS/Universit\' e Claude Bernard Lyon 1\\
43 Boulevard du 11 Novembre 1918,
F-69622 Villeurbanne Cedex, France\\
(kibler@ipnl.in2p3.fr)

\vspace*{.2cm} $^{4}$J. Heyrovsk\' y Institute of Physical
Chemistry, Academy of Sciences of the Czech Republic\\ Dolej\v
skova 3, CZ-182 23 Prague 8, Czech
Republic\\
(pracna@jh-inst.cas.cz)

\end{center}

\vspace*{-.2cm} \noindent \hrulefill

\vspace*{.1cm} \noindent {\bf Abstract}

\noindent A compact classification of the projective lines defined
over (commutative) rings (with unity) of all orders up to thirty-one is given.
There are altogether sixty-five different types
of them. For each type we introduce the total number of
points on the line, the number of points represented by
coordinates with at least one entry being a unit, the cardinality
of the neighbourhood of a generic point of the line as well as
those of the intersections between the neighbourhoods of two and
three mutually distant points, the number of `Jacobson' points per
a neighbourhood, the maximum number of pairwise distant points
and, finally, a list of representative/base rings. The
classification is presented in form of a table in order to see
readily not only the fine traits of the hierarchy, but also the
changes in the structure of the lines as one goes from one type to
the other. We hope this study will serve as an impetus to a search
for possible applications of these remarkable geometries in
physics, chemistry, biology and other natural sciences as well.
\\

%\vspace*{.15cm} \noindent {\bf MSC Codes:} 51C05 -- 81R99 -- 81Q99

%\vspace*{.05cm} \noindent {\bf PACS Numbers:} 02.10.Hh -- 02.40.Dr -- 03.65.Ca

%\vspace*{.05cm}

\vspace*{-.1cm} \noindent \hrulefill

\vspace*{.3cm}  \noindent
\section{Introduction}
Recently, projective lines defined over rings instead of fields
have become the subject of renewed interest not only for their own
sake \cite{bh1}--\cite{hav}, but also in view of their interesting
potential applications in (quantum) physics
\cite{sploc1}--\cite{sploc3}. It was the latter fact that
motivated our in-depth study of the structure of projective lines
over a large number of distinct commutative rings with unity of
order up to thirty-one
--- the study that yields, as far as we know, the first compact
classification of these beautiful geometric configurations. The key
element of our classification scheme is, of course, the neighbour
(or parallel) relation \cite{vk95},\cite{her}, which is a
geometrical concept intimately related with the structure of and
the connection between the maximal ideals of a ring under
consideration. Being simply an identity relation for fields, this
concept acquires a highly non-trivial character if the base ring
features three and/or more maximal ideals and, as we shall see in
detail, endows the corresponding line with a rich and quite
involved intrinsic structure.

\section{The Projective Line over a Ring with Unity}
Given a ring $R$  with unity (see, e.g., \cite{fr}--\cite{ra} and
also \cite{sploc1} or \cite{sploc3} for a brief recollection of
basic definitions, concepts and notations of ring theory) and
GL($2,R$), the general linear group of invertible two-by-two
matrices with entries in $R$, a pair ($\alpha, \beta$) $\in R^{2}$
is called {\it admissible} over $R$ if there exist $\gamma, \delta
\in R$ such that \cite{her}
\begin{equation}
\left(
\begin{array}{cc}
\alpha & \beta \\
\gamma & \delta \\
\end{array}
\right) \in {\rm GL}(2,R).
\end{equation}
The projective line over $R$, denoted as $PR(1)$, is defined as
the set of classes of ordered pairs $(\varrho \alpha, \varrho
\beta)$, where $\varrho$ is a unit and $(\alpha, \beta)$ is
admissible \cite{bh1},\cite{bh2},\cite{hav},\cite{her}. Such a
line carries two non-trivial, mutually complementary relations of
neighbour and distant. In particular, its two distinct points
$X$:=$(\varrho \alpha, \varrho \beta)$ and $Y$:=$(\varrho \gamma,
\varrho \delta)$ are called {\it neighbour} (or, {\it
parallel}) if
\begin{equation}
\left(
\begin{array}{cc}
\alpha & \beta \\
\gamma & \delta \\
\end{array}
\right) \notin {\rm GL}(2,R)
\end{equation}
and {\it distant} otherwise, i.\,e., if condition (1) is met. The
neighbour relation is reflexive (every point is obviously
neighbour to itself) and symmetric (i.\,e., if $X$ is neighbour to
$Y$ then $Y$ is neighbour to $X$ too), but, in general, not
transitive (i.\,e., $X$  being neighbour to $Y$ and $Y$ being
neighbour to $Z$ does not necessarily mean that $X$ is neighbour
to $Z$ --- see, e.\,g., \cite{hav},\cite{vk95},\cite{her}). Given
a point of $PR(1)$, the set of all neighbour points to it will be
called its {\it neighbourhood}.\footnote{To avoid any confusion,
the reader should be cautious that some authors (e.\,g.
\cite{bh1},\cite{hav}) use this term for the set of {\it distant}
points instead.} Obviously, if $R£$ is a field then `neighbour'
simply reduces to `identical' and `distant' to `different'. For a
{\it finite commutative} ring $R$, Eq.\,(1) reads
\begin{equation}
 \det \left(
\begin{array}{cc}
\alpha & \beta \\
\gamma & \delta \\
\end{array}
\right) \in R^{*},
\end{equation}
and Eq.\,(2) reduces to
\begin{equation}
\det \left(
\begin{array}{cc}
\alpha & \beta \\
\gamma & \delta \\
\end{array}
\right) \in R \backslash R^{*},
\end{equation}
where $R^{*}$ denotes the set of {\it units} (invertible elements)
and $R \backslash R^{*}$ stands for the set of {\it zero-divisors}
of $R$ (including the trivial zero divisor, 0).

To illustrate the concept of a ring line, we shall examine in
detail the structure of the projective line defined over the
direct product ring $R_{\diamondsuit} \equiv Z_{4}\otimes Z_{4}$,
with $Z_{4}$ being the ring of integers modulo 4, i.e., the set
$\{0, 1, 2, 3\}$ endowed with the addition and multiplication
properties as shown in Table 1.

\vspace*{-.4cm}
\begin{table}[h]
\begin{center}
\caption{Addition ({\it left}) and multiplication ({\it right}) in
$Z_4$.} \vspace*{0.3cm}
\begin{tabular}{||c|cccc||}
\hline \hline
$\oplus$ & 0 & 1 & 2 & 3 \\
\hline
0 & 0 & 1 & 2 & 3  \\
1 & 1 & 2 & 3 & 0  \\
2 & 2 & 3 & 0 & 1 \\
3 & 3 & 0 & 1 & 2 \\
\hline \hline
\end{tabular}~~~~~~~~~~~~~~~~~
\begin{tabular}{||c|cccc||}
\hline \hline
$\otimes$ & 0 & 1 & 2 & 3  \\
\hline
0 & 0 & 0 & 0 & 0  \\
1 & 0 & 1 & 2 & 3  \\
2 & 0 & 2 & 0 & 2  \\
3 & 0 & 3 & 2 & 1  \\
\hline \hline
\end{tabular}~.
\end{center}
\end{table}

\vspace*{-.4cm} \noindent The ring $R_{\diamondsuit}$ is, like
$Z_{4}$ itself, of characteristic four, and features the following
sixteen elements
\begin{eqnarray}
R_{\diamondsuit} = \{a \equiv [0,0], b \equiv [0,1], c \equiv
[0,2], d \equiv [0,3],
\nonumber \\
e \equiv [1,0], h \equiv [1,1], i \equiv [1,2], j \equiv [1,3],
\nonumber \\
f \equiv [2,0], k \equiv [2,1], l \equiv [2,2], m \equiv [2,3],
\nonumber \\
g \equiv [3,0], n \equiv [3,1], p \equiv [3,2], q \equiv [3,3] \}.
\end{eqnarray}
It contains two (proper) maximal ideals,
\begin{equation}
{\cal I}_{1} = \{ a, c, f, l, b, d, k, m \},
\end{equation}
\begin{equation}
{\cal I}_{2} = \{ a, c, f, l, e, g, i, p \},
\end{equation}
yielding a non-trivial Jacobson radical
\begin{equation}
{\cal J}= {\cal I}_{1}\cap {\cal I}_{2} = \{ a, c, f, l \},
\end{equation}
as it can readily be ascertained from its addition and
multiplication properties (Table 2).
\begin{table}[h]
\begin{center}
\caption{Addition ({\it top}) and multiplication ({\it bottom}) in
$R_{\diamondsuit}$.} \vspace*{0.2cm}
\begin{tabular}{||c|cccccccccccccccc||}
\hline \hline
$\oplus$ & $a$ & $b$ & $c$ & $d$ & $e$ & $f$ & $g$ & $h$ & $i$ & $j$ & $k$ & $l$ & $m$ & $n$ & $p$ & $q$\\
\hline
$a$ & $a$ & $b$ & $c$ & $d$ & $e$ & $f$ & $g$ & $h$ & $i$ & $j$ & $k$ & $l$ & $m$ & $n$ & $p$ & $q$ \\
$b$ & $b$ & $c$ & $d$ & $a$ & $h$ & $k$ & $n$ & $i$ & $j$ & $e$ & $l$ & $m$ & $f$ & $p$ & $q$ & $g$\\
$c$ & $c$ & $d$ & $a$ & $b$ & $i$ & $l$ & $p$ & $j$ & $e$ & $h$ & $m$ & $f$ & $k$ & $q$ & $g$ & $n$\\
$d$ & $d$ & $a$ & $b$ & $c$ & $j$ & $m$ & $q$ & $e$ & $h$ & $i$ & $f$ & $k$ & $l$ & $g$ & $n$ & $p$\\
$e$ & $e$ & $h$ & $i$ & $j$ & $f$ & $g$ & $a$ & $k$ & $l$ & $m$ & $n$ & $p$ & $q$ & $b$ & $c$ & $d$\\
$f$ & $f$ & $k$ & $l$ & $m$ & $g$ & $a$ & $e$ & $n$ & $p$ & $q$ & $b$ & $c$ & $d$ & $h$ & $i$ & $j$\\
$g$ & $g$ & $n$ & $p$ & $q$ & $a$ & $e$ & $f$ & $b$ & $c$ & $d$ & $h$ & $i$ & $j$ & $k$ & $l$ & $m$\\
$h$ & $h$ & $i$ & $j$ & $e$ & $k$ & $n$ & $b$ & $l$ & $m$ & $f$ & $p$ & $q$ & $g$ & $c$ & $d$ & $a$\\
$i$ & $i$ & $j$ & $e$ & $h$ & $l$ & $p$ & $c$ & $m$ & $f$ & $k$ & $q$ & $g$ & $n$ & $d$ & $a$ & $b$\\
$j$ & $j$ & $e$ & $h$ & $i$ & $m$ & $q$ & $d$ & $f$ & $k$ & $l$ & $g$ & $n$ & $p$ & $a$ & $b$ & $c$\\
$k$ & $k$ & $l$ & $m$ & $f$ & $n$ & $b$ & $h$ & $p$ & $q$ & $g$ & $c$ & $d$ & $a$ & $i$ & $j$ & $e$\\
$l$ & $l$ & $m$ & $f$ & $k$ & $p$ & $c$ & $i$ & $q$ & $g$ & $n$ & $d$ & $a$ & $b$ & $j$ & $e$ & $h$\\
$m$ & $m$ & $f$ & $k$ & $l$ & $q$ & $d$ & $j$ & $g$ & $n$ & $p$ & $a$ & $b$ & $c$ & $e$ & $h$ & $i$\\
$n$ & $n$ & $p$ & $q$ & $g$ & $b$ & $h$ & $k$ & $c$ & $d$ & $a$ & $i$ & $j$ & $e$ & $l$ & $m$ & $f$\\
$p$ & $p$ & $q$ & $g$ & $n$ & $c$ & $i$ & $l$ & $d$ & $a$ & $b$ & $j$ & $e$ & $h$ & $m$ & $f$ & $k$\\
$q$ & $q$ & $g$ & $n$ & $p$ & $d$ & $j$ & $m$ & $a$ & $b$ & $c$ & $e$ & $h$ & $i$ & $f$ & $k$ & $l$\\
\hline \hline
%\end{tabular}
&&&&&&&&&&&&&&&&\\
%\begin{tabular}{||c|cccccccccccccccc||}
\hline \hline
$\otimes$ &
      $a$ & $b$ & $c$ & $d$ & $e$ & $f$ & $g$ & $h$ & $i$ & $j$ & $k$ & $l$ & $m$ & $n$ & $p$ & $q$\\
\hline
$a$ & $a$ & $a$ & $a$ & $a$ & $a$ & $a$ & $a$ & $a$ & $a$ & $a$ & $a$ & $a$ & $a$ & $a$ & $a$ & $a$ \\
$b$ & $a$ & $b$ & $c$ & $d$ & $a$ & $a$ & $a$ & $b$ & $c$ & $d$ & $b$ & $c$ & $d$ & $b$ & $c$ & $d$\\
$c$ & $a$ & $c$ & $a$ & $c$ & $a$ & $a$ & $a$ & $c$ & $a$ & $c$ & $c$ & $a$ & $c$ & $c$ & $a$ & $c$\\
$d$ & $a$ & $d$ & $c$ & $b$ & $a$ & $a$ & $a$ & $d$ & $c$ & $b$ & $d$ & $c$ & $b$ & $d$ & $c$ & $b$\\
$e$ & $a$ & $a$ & $a$ & $a$ & $e$ & $f$ & $g$ & $e$ & $e$ & $e$ & $f$ & $f$ & $f$ & $g$ & $g$ & $g$\\
$f$ & $a$ & $a$ & $a$ & $a$ & $f$ & $a$ & $f$ & $f$ & $f$ & $f$ & $a$ & $a$ & $a$ & $f$ & $f$ & $f$\\
$g$ & $a$ & $a$ & $a$ & $a$ & $g$ & $f$ & $e$ & $g$ & $g$ & $g$ & $f$ & $f$ & $f$ & $e$ & $e$ & $e$\\
$h$ & $a$ & $b$ & $c$ & $d$ & $e$ & $f$ & $g$ & $h$ & $i$ & $j$ & $k$ & $l$ & $m$ & $n$ & $p$ & $q$\\
$i$ & $a$ & $c$ & $a$ & $c$ & $e$ & $f$ & $g$ & $i$ & $e$ & $i$ & $l$ & $f$ & $l$ & $p$ & $g$ & $p$\\
$j$ & $a$ & $d$ & $c$ & $b$ & $e$ & $f$ & $g$ & $j$ & $i$ & $h$ & $m$ & $l$ & $k$ & $q$ & $p$ & $n$\\
$k$ & $a$ & $b$ & $c$ & $d$ & $f$ & $a$ & $f$ & $k$ & $l$ & $m$ & $b$ & $c$ & $d$ & $k$ & $l$ & $m$\\
$l$ & $a$ & $c$ & $a$ & $c$ & $f$ & $a$ & $f$ & $l$ & $f$ & $l$ & $c$ & $a$ & $c$ & $l$ & $f$ & $l$\\
$m$ & $a$ & $d$ & $c$ & $b$ & $f$ & $a$ & $f$ & $m$ & $l$ & $k$ & $d$ & $c$ & $b$ & $m$ & $l$ & $k$\\
$n$ & $a$ & $b$ & $c$ & $d$ & $g$ & $f$ & $e$ & $n$ & $p$ & $q$ & $k$ & $l$ & $m$ & $h$ & $i$ & $j$\\
$p$ & $a$ & $c$ & $a$ & $c$ & $g$ & $f$ & $e$ & $p$ & $g$ & $p$ & $l$ & $f$ & $l$ & $i$ & $e$ & $i$\\
$q$ & $a$ & $d$ & $c$ & $b$ & $g$ & $f$ & $e$ & $q$ & $p$ & $n$ & $m$ & $l$ & $k$ & $j$ & $i$ & $h$\\
\hline \hline
\end{tabular}~.
\end{center}
\end{table}

\noindent From these tables it also follows that $a$ and $h$ are,
respectively, the addition and multiplication identities (`0' and
`1') of the ring  and that
\begin{equation}
R^{*}_{\diamondsuit} = \{h \equiv 1, j, n, q\}
\end{equation}
and
\begin{equation}
R_{\diamondsuit} \backslash R^{*}_{\diamondsuit} = \{a \equiv 0,
b, c, d, e, f, g, i, k, l, m, p \}.
\end{equation}
Now we can employ the admissibility condition and find out that
the projective line $PR_{\diamondsuit}(1)$ contains altogether
thirty-six points, which can be partitioned into two distinct
groups. The first group (`type I points') consists of the points
represented by coordinates where at least one entry is a unit and
there are the following twenty-eight points here
\begin{eqnarray}
&&(1,0), (1,b),(1,c),(1,d),(1,e),(1,f),(1,g), (1,i), (1,k),(1,l),
(1,m), (1,p),
\nonumber \\
&& (0,1), (b,1),(c,1),(d,1),(e,1),(f,1),(g,1), (i,1), (k,1),(l,1),
(m,1), (p,1), \nonumber \\
&& (1,1), (1,j), (1,n), (1,q);
\end{eqnarray}
it is easy to verify that for any finite commutative ring this
number is always equal to the sum of the total elements of the
ring and the number of its zero-divisors. The other group (`type
II points') is composed of those points whose representing
coordinates are both zero-divisors and we find the following eight
points ranked here
\begin{equation}
(e,b), (e,k), (i,b), (i,k), (b,e), (k,e), (b,i), (k,i);
\end{equation}
these points, in general, exist only if the ring has two or more maximal ideals
and their number depends on the properties and interconnection
between these ideals. To reveal all the subtleties of the
structure of the line, one has to make use of the
neighbour/distant relation. The reasoning here is, without any
loss of generality, much facilitated by considering three
distinguished points of the line, viz.
$
U:=(1, 0),~~ V:=(0, 1),~~ W:=(1, 1),
$
which are obviously pairwise distant and have the following
neighbourhoods
\begin{eqnarray}
U:&&(1,b), \underline{(1,c)}, (1,d), (1,e),
\underline{(1,f)},(1,g), (1,i), (1,k),\underline{(1,l)}, (1,m),
(1,p), \nonumber \\
&&(e,b), (e,k), (i,b), (i,k), (b,e), (k,e), (b,i), (k,i)
\end{eqnarray}
\begin{eqnarray}
V:&&(b,1), \underline{(c,1)}, (d,1), (e,1),
\underline{(f,1)},(g,1), (i,1), (k,1),\underline{(l,1)}, (m,1),
(p,1), \nonumber \\
&&(e,b), (e,k), (i,b), (i,k), (b,e), (k,e), (b,i), (k,i)
\end{eqnarray}
\begin{eqnarray}
W: &&(1,b),(1,d), (1,e), (1,g), (1,i), (1,k), (1,m),
(1,p), \nonumber \\
&& (b,1), (d,1), (e,1),(g,1), (i,1), (k,1),(m,1),
(p,1), \nonumber \\
&& \underline{(1,j)}, \underline{(1,n)}, \underline{(1,q)}.
\end{eqnarray}
We see that each neighbourhood features nineteen points and has
three `Jacobson' points (underlined), i.e., the points unique to
the particular neighbourhood, that the neighbourhoods pairwise
overlap in eight points and have no common element if considered
altogether; moreover, one easily checks that there exists no point
of the line that would be simultaneously distant to all the three
distinguished points. Now, as the coordinate system on this line
can {\it always} be chosen in such a way that the coordinates of
{\it any} three mutually distant points are made identical to
those of $U$, $V$ and $W$, we can conclude that the neighbourhood
of any point of the line features nineteen distinct points, the
neighbourhoods of any two distant points share eight points, the
neighbourhoods of any three mutually distant points have no point
in common and three is the maximum number of mutually distant
points; a nice `conic' representation of the line exhibiting all
these properties is given in Fig.\,1.
\begin{figure}[h]
\centerline{\includegraphics[width=7.0truecm,clip=]{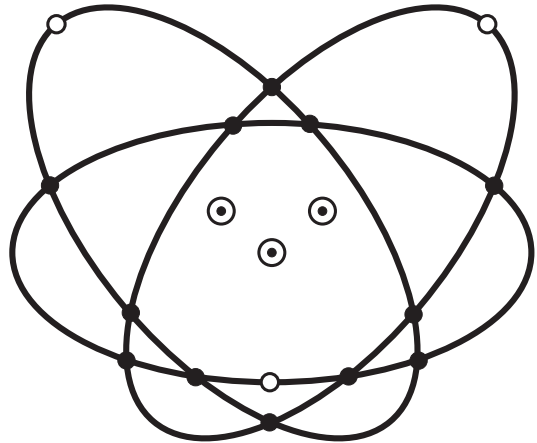}}
\caption{A schematic sketch of the structure of the projective
line $PR_{\diamondsuit}(1)$. The three double circles stand for
three pairwise distant points, whilst the remaining points of the
line are all located on the neighbourhoods of these three points
(three sets of points located on three different ellipses centered
at the points in question). Every bullet represents {\it two}
distinct points of the line, while each of the three small circles
represents {\it three} `Jacobson' points.}
\end{figure}

\section{Classifying Lines over Rings of Small Orders}

Following the same strategy as in the previous
section, we have examined the properties of the projective lines
over commutative rings with unity of all orders up to
thirty-one. The most relevant results/findings of our study are in a
succinct, and in a number of respects also unique, form displayed
in Table 3; note that the line analysed in the preceding section is
of type 16/12. The following abbreviations are used for the
number of points: `Tot' -- all the points of the line; `TpI' --
points of type I; `1N' -- in the neighbourhood of a generic point
(the point itself exclusive); `$\cap$2N'-- common to the
neighbourhoods of two distant points; `$\cap$3N'-- shared by the
neighbourhoods of three pairwise distant points; `Jcb' -- having
`Jacobson' property; `MD' -- the maximum set of pairwise distant
ones. The type of a line is given in the form A/B, where A is the
total number of elements of the associated ring and B is the
number of its zero-divisors. In all cases except for 8/4, 16/8,
16/12, 24/16 and 27/9 the table gives all the representative rings; due to a lack
of space here, a full list of rings for each of the five
specified types is given in a separate table (Table 4), using the
notation of \cite{cn1} (orders $\leq$ 16) and \cite{cn2} (orders $>$16).

\vspace*{1.cm}
\begin{table}[hb]
\begin{center}
\caption{The basic types of a projective line over small commutative rings with unity. For the types denoted
by bold-facing
there also exist ``non-commutative" projective lines of the particular orders \cite{spp}.}
\vspace*{0.8cm}
{\begin{tabular}{|||l|r|r|r|r|r|c|c|l|||} \hline \hline \hline
\multicolumn{1}{|||c|}{} & \multicolumn{7}{|c|}{} &  \multicolumn{1}{|c|||}{}\\
%\cline{2-13}
\multicolumn{1}{|||c|}{Line} & \multicolumn{7}{|c|}{Cardinalities
of Points} &
\multicolumn{1}{|c|||}{Representative} \\
\multicolumn{1}{|||c|}{Type} & \multicolumn{7}{|c|}{} &  \multicolumn{1}{|c|||}{Rings} \\
\cline{2-8}
 &Tot& TpI & 1N & $\cap$2N & $\cap$3N & Jcb & MD &   \\
\hline
\hline
31/1 & 32 & 32 & 0 & 0 & 0 & 0 & 32 & $GF(31)$  \\
\hline \hline
30/22 & 72 & 52 & 41 & 20 & 6 & 7 & 3 & $GF(2)\otimes GF(3)\otimes GF(5) $  \\
\hline \hline
29/1 & 30 & 30 & 0 & 0 & 0 & 0 & 30 & $GF(29)$  \\
\hline \hline
28/10 & 40 & 38 & 11 & 2 & 0 & 3 & 5 & $GF(4)\otimes GF(7)$  \\
\hline
28/16 & 48 & 44 & 19 & 4 & 0 & 11 & 3 & $GF(7)\otimes
[Z_{4}$ or $GF(2)[x]/\langle x^{2}\rangle$]\\
28/22 & 72 & 50 & 43 & 22 & 6 & 5 & 3 & $GF(2)\otimes GF(2)\otimes GF(7)$  \\
\hline \hline
27/1 & 28 & 28 & 0 & 0 & 0 & 0 & 28 & $GF(27)$  \\
\hline
27/9 & 36 & 36 & 8 & 0 & 0 & 8 & 4 & $Z_{27}, GF(3)[x]/\langle x^{3}\rangle, \ldots$  \\
27/11 & 40 & 38 & 12 & 2 & 0 & 6 & 4 & $GF(3)\otimes GF(9)$  \\
{\bf 27/15} & 48 & 42 & 20 & 6 & 0 & 2 & 4 & $GF(3)\otimes [Z_{9}$ or $GF(3)[x]/\langle x^{2}\rangle$]  \\
27/19 & 64 & 46 & 36 & 18 & 6 & 0 & 4 & $GF(3)\otimes GF(3)\otimes GF(3)$  \\
\hline \hline
26/14 & 42 & 40 & 15 & 2 & 0 & 11 & 3 & $GF(2)\otimes GF(13)$  \\
\hline \hline
25/1 & 26 & 26 & 0 & 0 & 0 & 0 & 26 & $GF(25)$  \\
\hline
25/5 & 30 & 30 & 4 & 0 & 0 & 4 & 6 & $Z_{25}, GF(5)[x]/\langle x^{2}\rangle$  \\
25/9 & 36 & 34 & 10 & 2 & 0 & 0 & 6 & $GF(5)\otimes GF(5)$  \\
\hline \hline
24/10 & 36 & 34 & 11 & 2 & 0 & 5 & 4 & $GF(3)\otimes GF(8)$  \\
\hline
24/16 & 48 & 40 & 23 & 8 & 0 & 7 & 3 & $GF(3)\otimes Z_{8}, GF(3)\otimes GF(2)[x]/\langle x^{3}\rangle, \ldots $  \\
24/18 & 60 & 42 & 35 & 18 & 6 & 5 & 3 & $GF(2)\otimes GF(3)\otimes GF(4) $  \\
{\bf 24/20} & 72 & 44 & 47 & 28 & 12 & 3 & 3 & $GF(2)\otimes GF(3)\otimes [Z_{4}$ or $GF(2)[x]/\langle x^{2}\rangle$]   \\
24/22 & 108 & 46 & 83 & 62 & 42 & 1 & 3 & $GF(2)\otimes GF(2)\otimes GF(2)\otimes GF(3) $  \\
\hline
\hline
23/1 & 24 & 24 & 0 & 0 & 0 & 0 & 24 & $GF(23)$  \\
\hline \hline
 22/12 & 36 & 34 & 13 & 2 & 0 & 9 & 3 & $GF(2)\otimes GF(11)$  \\
\hline \hline
 21/9 & 32 & 30 & 10 & 2 & 0 & 4 & 4 & $GF(3)\otimes GF(7)$  \\
\hline \hline
20/8 & 30 & 28 & 9 & 2 & 0 & 1 & 5 & $GF(5)\otimes GF(4)$   \\
\hline
20/12 & 36 & 32 & 15 & 4 & 0 & 7 & 3 & $GF(5)\otimes [Z_{4}$ or $GF(2)[x]/\langle x^{2}\rangle$]  \\
%\hline
20/16 & 54 & 36 & 33 & 18 & 6 & 3 & 3 & $GF(5)\otimes GF(2)\otimes GF(2)$  \\
\hline
\hline
 19/1 & 20 & 20 & 0 & 0 & 0 & 0 & 20 & $GF(19)$  \\
\hline \hline
18/10 & 30 & 28 & 11 & 2 & 0 & 7 & 3 &  $GF(2)\otimes GF(9)$ \\
18/12 & 36 & 30 & 17 & 6 & 0 & 5 & 3 &  $GF(2)\otimes [Z_{9}$ or $GF(3)[x]/\langle x^{2}\rangle$] \\
18/14 & 48 & 32 & 29 & 16 & 6 & 3 & 3 &  $GF(2)\otimes GF(3)\otimes GF(3)$ \\
\hline
\hline
17/1 & 18 & 18 & 0 & 0 & 0 & 0 & 18 & $GF(17)$  \\
\hline
\hline \hline
\end{tabular}}
\end{center}
\end{table}
\clearpage
\addtocounter{table}{-1}

\begin{table}[ht]
\begin{center}
\caption{(Continued)}
\vspace*{0.3cm}
{\begin{tabular}{|||l|r|r|r|r|r|c|c|l|||} \hline \hline \hline
\multicolumn{1}{|||c|}{} & \multicolumn{7}{|c|}{} &  \multicolumn{1}{|c|||}{}\\
%\cline{2-13}
\multicolumn{1}{|||c|}{Line} & \multicolumn{7}{|c|}{Cardinalities
of Points} &
\multicolumn{1}{|c|||}{Representative} \\
\multicolumn{1}{|||c|}{Type} & \multicolumn{7}{|c|}{} &  \multicolumn{1}{|c|||}{Rings} \\
\cline{2-8}
 &Tot& TpI & 1N & $\cap$2N & $\cap$3N & Jcb & MD &   \\
\hline
\hline
16/1 & 17 & 17 & 0 & 0 & 0 & 0 & 17 & $GF(16)$  \\
\hline
{\bf 16/4} & 20 & 20 & 3 & 0 & 0 & 3 & 5 & $Z_{4}[x]/\langle x^{2} + x + 1\rangle$, $GF(4)[x]/\langle x^{2}\rangle$  \\
%\hline
16/7 & 25 & 23 & 8 & 2 & 0 & 0 & 5 & $GF(4)\otimes GF(4)$ \\
\hline
{\bf 16/8} & 24 & 24 & 7 & 0 & 0 & 7 & 3 & $Z_{16}$,
$Z_{4}[x]/\langle
x^{2}\rangle$, $GF(2)[x]/\langle x^{4}\rangle$, \ldots  \\
%\hline
16/9 & 27 & 25 & 10 & 2 & 0 & 6 & 3 & $GF(2)\otimes GF(8)$
\\
%\hline
{\bf 16/10} & 30 & 26 & 13 & 4 & 0 & 5 & 3 & $GF(4)\otimes
[Z_{4}$ or $GF(2)[x]/\langle x^{2}\rangle$]\\
%\hline
{\bf 16/12} & 36 & 28 & 19 & 8 & 0 & 3 & 3 &  $Z_{4} \otimes
Z_{4}$, $GF(2)\otimes Z_{8}$, \ldots \\
%\hline
16/13 & 45 & 29 & 28 & 16 & 6 & 2 & 3 &  $GF(2)\otimes GF(2)\otimes GF(4)$\\
%\hline
{\bf 16/14} & 54 & 30 & 37 & 24 & 12 & 1 & 3 &  $GF(2)\otimes GF(2)\otimes [Z_{4} ~{\rm or}~ GF(2)[x]/\langle x^{2}\rangle$] \\
%\hline
16/15 & 81 & 31 & 64 & 50 & 36 & 0 & 3 &  $GF(2)\otimes GF(2)\otimes GF(2)\otimes GF(2)$ \\
\hline
\hline
15/7 & 24 & 22 & 8 & 2 & 0 & 2 & 4 &  $GF(3)\otimes GF(5)$ \\
\hline
\hline
14/8 & 24 & 22 & 9 & 2 & 0 & 5 & 3 &  $GF(2)\otimes GF(7)$ \\
\hline
\hline
13/1 & 14 & 14 & 0 & 0 & 0 & 0 & 14 & $GF(13)$  \\
\hline
\hline
12/6 & 20 & 18 & 7 & 2 & 0 & 1 & 4 &  $GF(3)\otimes GF(4)$ \\
\hline
12/8 & 24 & 20 & 11 & 4 & 0 & 3 & 3 &  $GF(3)\otimes [Z_{4} ~{\rm or}~ GF(2)[x]/\langle x^{2}\rangle$] \\
12/10 & 36 & 22 & 23 & 14 & 6 & 1 & 3 &  $GF(2) \otimes GF(2) \otimes GF(3)$ \\
\hline
\hline
11/1 & 12 & 12 & 0 & 0 & 0 & 0 & 12 & $GF(11)$  \\
\hline
\hline
10/6 & 18 & 16 & 7 & 2 & 0 & 3 & 3 & $GF(2)\otimes GF(5)$  \\
\hline
\hline
\hspace*{0.2cm}9/1 & 10 & 10 & 0 & 0 & 0 & 0 & 10 &  $GF(9)$ \\
\hline
\hspace*{0.2cm}9/3 & 12 & 12 & 2 & 0 & 0 & 2 & 4 &  $Z_{9}$, $GF(3)[x]/\langle x^{2}\rangle$ \\
\hspace*{0.2cm}9/5 & 16 & 14 & 6 & 2 & 0 & 0 & 4 &  $GF(3)\otimes GF(3)$ \\
\hline
\hline
\hspace*{0.2cm}8/1 & 9 & 9 & 0 & 0 & 0 & 0 & 9 & $GF(8)$  \\
\hline \hspace*{0.2cm}8/4 & 12 & 12 & 3 & 0 & 0 & 3 & 3 & $Z_{8}, GF(2)[x]/\langle x^{3}\rangle, \ldots$  \\
%\hline
\hspace*{0.2cm}8/5 & 15 & 13 & 6 & 2 & 0 & 2 & 3 & $GF(2)\otimes GF(4)$  \\
%\hline
\hspace*{0.2cm}{\bf 8/6} & 18 & 14 & 9 & 4 & 0 & 1 & 3 & $GF(2)\otimes [Z_{4}$ or $GF(2)[x]/\langle x^{2}\rangle$]  \\
%\hline
\hspace*{0.2cm}8/7 & 27 & 15 & 18 & 12 & 6 & 0 & 3 & $GF(2)\otimes GF(2)\otimes GF(2)$  \\
\hline
\hline
\hspace*{0.2cm}7/1 & 8 & 8 & 0 & 0 & 0 & 0 & 8 & $GF(7)$  \\
\hline
\hline
\hspace*{0.2cm}6/4 & 12 & 10 & 5 & 2 & 0 & 1 & 3 & $GF(2)\otimes GF(3)$  \\
\hline
\hline
\hspace*{0.2cm}5/1 & 6 & 6 & 0 & 0 & 0 & 0 & 6 & $GF(5)$  \\
\hline
\hline
\hspace*{0.2cm}4/1 & 5 & 5 & 0 & 0 & 0 & 0 & 5 & $GF(4)$  \\
\hline \hspace*{0.2cm}4/2 & 6 & 6 & 1 & 0 & 0 & 1 & 3 & $Z_{4}$, $GF(2)[x]/\langle x^{2}\rangle$  \\
%\hline
\hspace*{0.2cm}4/3 & 9 & 7 & 4 & 2 & 0 & 0 & 3 & $GF(2)\otimes GF(2)$  \\
\hline
\hline
\hspace*{0.2cm}3/1 & 4 & 4 & 0 & 0 & 0 & 0 & 4 & $GF(3)$  \\
\hline
\hline
\hspace*{0.2cm}2/1 & 3 & 3 & 0 & 0 & 0 & 0 & 3 & $GF(2)$  \\
\hline \hline \hline\end{tabular}}
\end{center}
\end{table}

\begin{table}[htb]
\begin{center}
\caption{A comprehensive list of all commutative rings with unity
which generate the projective lines of the type 8/4, 16/8,
16/12, 24/16 and 27/9.} \vspace*{0.8cm}
{\begin{tabular}{||l|l||} \hline \hline  Type of Line &
Representative
Rings \\
\hline \hline \hspace*{0.2cm}8/4 & 1.4, 2.19, 2.20, 3.14, 3.19 \\
\hline 16/8 & 1.5, 2.28, 2.29, 3.40, 3.42, 3.61, 3.64,
4.118, 4.119, 4.165, 4.166, 4.168, 4.170, \\
& 5.99, 5.102, 5.109, 5.110, 5.111  \\
\hline 16/12 & 2.19, 3.49, 4.81, 4.82, 4.167, 5.107, 5.114, 5.115
\\
\hline 24/16 & 1.8, 2.39, 2.40, 3.28, 3.38
\\
\hline 27/9 & 1.4, 2.17, 2.24, 2.25, 3.14, 3.23
\\
\hline \hline
\end{tabular}}
\end{center}
\end{table}

\noindent
From Table 3 one can discern a number
of interesting features of the structure of this hierarchy. The
most marked one is a general increase of the total number of
points, those of type II and the population of neighbourhoods with
the increasing number of zero-divisors of the base ring. This is
accompanied by strengthening of the `coupling,' or `entanglement,'
between the neighbourhoods of mutually distant points, which is
embodied in gradually increasing overlaps between the
neighbourhoods of first two and then three such points; note, for
example, that for lines of 16/15 type, the number of points in the
intersection of the neighbourhoods of three pairwise distant
points is greater than that of type I. One further observes that
various types of a projective line form natural groups (Gr)
differing from each other in the total number of elements of the
associated base ring. Each such group can further be divided into
classes (Cl) of the same maximum number of mutually distant
points; thus, for example, Gr-16 consists of three classes Cl-17,
Cl-5 and Cl-3, the last-mentioned being the most populated.
Classes generated by the Galois fields, $GF(q)$, contain just a
single entry and are regarded as trivial. A further subdivision
within a given class is into cells (Ce) according as the
neighbourhoods of mutually distant points are disjoint
(Ce$\cap$1), overlaps for pairs (Ce$\cap$2) and/or triples
(Ce$\cap$3) of them; so, for example, the class Cl-3 of Gr-16 is
seen to comprise all the three kinds of a cell, of the
cardinalities one, three and three, respectively, whilst the same
class of Gr-2 features just a single cell, Ce$\cap$1. A point that
deserves particular attention here is a considerable change in the
structure of the line under the transition between the classes
within a given group. This change is likely to be found ever more
pronounced with the increasing order of the ring and for lines
within our scope it is best visible in Gr-16, on the boundary
between its Cl-5 and Cl-3; one sees that drop, when moving from
16/7 to 16/8, in the maximum number of pairwise distant and the
total number of points, as well as in the cardinality of a
neighbourhood, is accompanied by reappearance of `Jacobson' points
and disappearance of the cell of type Ce$\cap$2. Another
noteworthy property is gradual decrease in the number of
`Jacobson' points within any (non-trivial) class. It is also to be
noted that the neighbourhoods of mutually distant points are
disjoint on the lines defined over local rings (types 4/2, 8/4, 9/3,
16/4, 16/8, 25/5 and 27/9), a fact which entails the transitivity of the
neighbour relation.\\

\section{Conclusion}
All projective lines defined over commutative rings with unity
of orders two to thirty-one have been classified. We
have found altogether sixty-five different types of them. Each type is characterized
by the following string of parameters: the total number of points
on the line, the number of points represented by coordinates with
at least one entry being a unit, the cardinality of the
neighbourhood of a generic point of the line as well as those of
the intersections between the neighbourhoods of two and three
mutually distant points, the number of `Jacobson' points per a
neighbourhood, the maximum number of pairwise distant points and,
finally, by a list of representative rings. The exposition of the
ideas and the classification itself are presented in the way to
stir the interest of physicists, chemists, biologists and scholars
of other natural sciences to look for possible applications of
these abstract finite geometries in their domains of research. As
per physics, a couple of the above-introduced types of a
projective ring line, namely the 4/2 and 8/6 ones, have already
been successfully employed to account for some subtleties of the
structure of {\it two-qubit} systems \cite{sploc1},\cite{sploc2};
our most recent study \cite{psk} indicates that it is also the
line of type 16/15 whose structure is relevant for these quantum
systems, with type II points (and, so, zero-divisors) playing a
particular role here. The standard model of elementary particles
and their interactions is another domain where the combinatorics
of finite ring lines is likely to find its proper setting; in this
respect, it should be emphasized that the line of type 4/2 could
be of interest for the classification of the six quarks and six
leptons.

\vspace*{.7cm} \noindent \Large {\bf Acknowledgements}
\normalsize

\vspace*{.3cm} \noindent This work was partially supported by the
Science and Technology Assistance Agency under the contract $\#$
APVT--51--012704, the VEGA project $\#$ 2/6070/26 (both from
Slovak Republic), the trans-national ECO-NET project $\#$
12651NJ ``Geometries Over Finite Rings and the Properties of
Mutually Unbiased Bases" (France) and by the project 1ET400400410 of
the Academy of Sciences of the Czech Republic.


\vspace*{-.1cm}
\begin{thebibliography}{10}
\vspace*{-.20cm}
\bibitem{bh1}
Blunck A, Havlicek H. Projective representations I: Projective
lines over a ring.  Abh Math Sem Univ Hamburg 2000;70:287--99.
\vspace*{-.25cm}
\bibitem{bh2}
Blunck A, Havlicek H. Radical parallelism on projective lines and
non-linear models of affine spaces. Mathematica Pannonica
2003;14:113--27. \vspace*{-.25cm}
\bibitem{bh3}
Blunck A, Havlicek H. On distant-isomorphisms of projective lines.
Aequationes Mathematicae 2005;69:146--163. \vspace*{-.25cm}
\bibitem{hav}
Havlicek H. Divisible designs, Laguerre geometry, and beyond.
Quaderni del Seminario Matematico di Brescia 2006;11:1--63. A
preprint available from
$\langle$http://www.geometrie.tuwien.ac.at/havlicek/dd-laguerre.pdf$\rangle$.
\vspace*{-.25cm}
\bibitem{sploc1}
Saniga M, Planat M. The projective line over the finite quotient
ring GF(2)[$x$]/$\langle x^3-x \rangle$ and quantum entanglement
I. Theoretical background. Theoretical and Mathematical Physics 2006; submitted. Available from
$\langle$quant-ph/0603051$\rangle$. \vspace*{-.25cm}
\bibitem{sploc2}
Saniga M, Planat M, Minarovjech M. The projective line over the
finite quotient ring \linebreak GF(2)[$x$]/$\langle x^3-x \rangle$
and quantum entanglement II. The Mermin ``Magic''
Square/Pentagram. Theoretical and Mathematical Physics 2006; submitted. Available from
$\langle$quant-ph/0603206$\rangle$. \vspace*{-.25cm}
\bibitem{sploc3}
Saniga M, Planat, M. On the fine structure of the projective line
over GF(2) $\otimes$ GF(2) $\otimes$ GF(2). Indagationes Mathematicae 2006; submitted.
Available from
$\langle$math.AG/0604307$\rangle$. \vspace*{-.25cm}
\bibitem{vk95}
Veldkamp FD. Geometry over rings. In: Buekenhout F, editor.
Handbook of incidence geometry.
        Amsterdam: Elsevier; 1995. p.\,1033--84.
\vspace*{-.25cm}
\bibitem{her} Herzer A. Chain geometries. In:
Buekenhout F, editor. Handbook of incidence geometry.
        Amsterdam: Elsevier; 1995. p.\,781--842.
        \vspace*{-.25cm}
\bibitem{fr}
Fraleigh JB. A first course in abstract algebra (5th edition).
Reading (MA): Addison-Wesley; 1994. p.\,273--362. \vspace*{-.25cm}
\bibitem{mcd}
McDonald BR. Finite rings with identity. New York: Marcel Dekker;
1974. \vspace*{-.25cm}
\bibitem{ra}
Raghavendran R. Finite associative rings. Comp Mathematica
1969;21:195--229.
\vspace*{-.25cm}
\bibitem{cn1}
N\" obauer C. The Book of the Rings --- Part I. 2000; Available
from
$\langle$http://www.algebra.uni-linz.ac.at/$\widetilde{~~}$noebsi/pub/rings.ps$\rangle$.
\vspace*{-.25cm}
\bibitem{cn2}
N\" obauer C. The Book of the Rings --- Part II. 2000; Available
from
$\langle$http://www.algebra.uni-linz.ac.at/$\widetilde{~~}$noebsi/pub/ringsII.ps$\rangle$.
\vspace*{-.25cm}
\bibitem{spp}
Saniga M, Planat M, Pracna P. A classification of the projective lines over small rings II. Non-commutative case.
Available from
$\langle$math.AG/0606500$\rangle$.
\vspace*{-.25cm}
\bibitem{psk}
Planat M, Saniga M, Kibler MR. Quantum entanglement and finite ring geometry. SIGMA 2006; accepted.
Available from
$\langle$quant-ph/0605239$\rangle$.
\end{thebibliography}
\end{document}